\documentclass[onecolumn,10pt]{article}
\usepackage{cite}
\usepackage{graphicx}
\usepackage[cmex10]{amsmath}
\usepackage[tight,footnotesize]{subfigure}
\usepackage{amsfonts,latexsym,amssymb,amsmath,amsthm}

\begin{document}
\title{Finding the Dominant Roots of a Time Delay System without Using the Principal Branch of the Lambert W Function}
\author{Rudy Cepeda-Gomez\\
\textit{a+i Engineering, Berlin, Germany}\\
\texttt{rudy.cepeda-gomez@a--i.com}}%
\maketitle
\begin{abstract}
This brief note complements some results regarding a recently developed technique for the stability analysis of linear time-invariant, time delay systems using the matrix Lambert W function. By means of a numeric example, it is shown that there are cases for which the dominant roots of the system can be found without using the principal branch of this multi-valued function, contradicting the main proposition of the methodology.
\end{abstract}
\section{Introduction}
Consider a $n$-th order Linear Time Invariant-Time Delay System (\textsc{lti-tds}), represented by a Delay-Differential Equation (\textsc{dde}) of the form
\begin{equation}
\dot{x}\left(t\right)=Ax\left(t\right)+Bx\left(t-\tau\right),
\label{eq:lti-tds}
\end{equation}
with $x\in\mathbb{R}^n$, $A,B\in\mathbb{R}^{n\times n}$ and $\tau>0$. A framework for analyzing \textsc{ddes} like \eqref{eq:lti-tds} based on the Lambert W function has been developed in the past decade \cite{Asl2003,Yi2009,Yi2010d}, expanding earlier works like \cite{Wright1959}. The main idea of this methodology is to express the solution of a \textsc{dde} as the sum of a series of infinitely many exponential functions. The characteristic roots of the system are found analytically in terms of the Lambert W function. While the problem remains infinite dimensional, a one-to-one correspondence between the characteristic roots of the system and the branches of this multi-valued function is assumed. The stability question is then solved by earmarking the dominant characteristic roots of the system with the branches of the Lambert W function, such that only a few branches have to be considered to determine whether a solution is stable or not.

Although the basic foundation of this methodology, i.e., the assumption that the principal branch of the Lambert W function defines the stability of the system, was proven for first order systems (when $n=1$ in \eqref{eq:lti-tds}) in \cite{Shinozaki2006} and \cite{Asl2003} independently, for higher order systems this is not the case. It was shown in \cite{AUT} that under mild conditions it is possible to find the complete spectrum of a second order system using only two branches of the matrix Lambert W function, namely $k=0$ and $k=-1$.

Nevertheless, the example presented in \cite{AUT} shows that the 22 dominant roots are found using the principal branch. This may prompt the reader to believe that there is still some of the correspondence observed for first order systems. In this technical note, we present an example in which \emph{all roots of the system can be found without using the principal branch of the matrix Lambert W function}. All the roots in this case can be found using the branch corresponding to $k=-1$.

To avoid unnecessary repetition, this note does not describe the steps of the stability analysis methodology based on the matrix Lambert W function, which can be found in references such as \cite{Yi2009,Yi2006a,Yi2006b}, nor the foundation and the steps used to create the counterexamples, which are presented in \cite{AUT}. For a discussion on the definition and properties of the Lambert W function, the reader is referred to \cite{Corless1996}.
\section{Finding the Spectrum Without using The Principal Branch}
\begin{figure}[tb]
\centering
\includegraphics[scale=1]{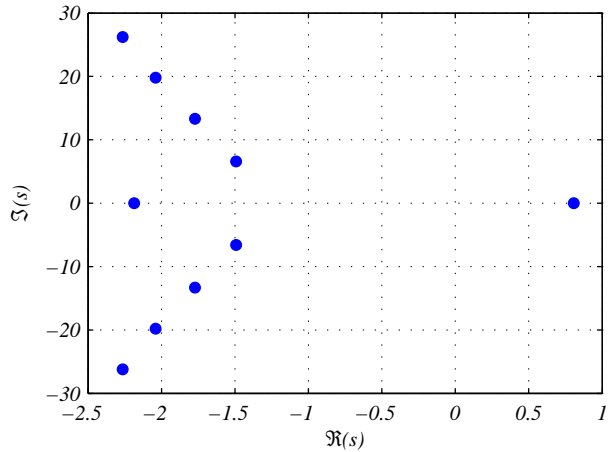}
\caption{Dominant roots of the system under study.}
\label{fig:rts}
\end{figure}

For a system described by \eqref{eq:lti-tds} with
\begin{equation}
A=\left[\begin{array}{rr}0&1\\-5&10\end{array}\right]\quad B=\left[\begin{array}{rr}0&0\\-3&-3\end{array}\right],
\end{equation}
and $\tau=1$, the QPmR algorithm \cite{Vhylidal2009} finds the 10 dominant roots observed in Fig.~\ref{fig:rts}. Following the reverse engineering approach described in \cite{AUT}, the two real roots of the system, $\lambda_1=0.8070$ and $\lambda_2=-2.1854$, are combined to create the matrix
\begin{equation}
S=\left[\begin{array}{cc}0&1\\1.7636&-1.3784\end{array}\right],
\end{equation}
which in turns generates
\begin{equation}
W_k(M_k)=\tau\left(S-A\right)=\left[\begin{array}{cc}0&0\\6.7636&-11.3784\end{array}\right].
\label{eq:wkm}
\end{equation}
As discussed in \cite{AUT}, since $W(m_{22})\in\left(-\infty,\,-1\right)$, which is the range of the branch $k=-1$ of the Lambert W function, there is an $M$ matrix such that \eqref{eq:wkm} is satisfied for $k=-1$. That matrix is
\begin{equation}
M_{-1}=\left[\begin{array}{cc}0&0\\0.0774&-0.1302\end{array}\right].
\end{equation}
and from it, a $Q_{-1}$ matrix like
\begin{equation}
Q_{-1}=\left[\begin{array}{rr}2&1\\-2&-1\end{array}\right].
\end{equation}
can be used as a starting value in the \emph{LambertDDE} Matlab toolbox \cite{Yi2012a} to obtain these roots as a solution with $k=-1$. Notice that this pair of roots includes the dominant root of the system, which according to the methodology under study should have been found using $k=0$.

If the second dominant pair of roots of the system, i.e. $s=-1.4928\pm j6.6027$, is used to create a matrix like
\begin{equation}
S=\left[\begin{array}{cc}0&1\\-45.8241&-2.9855\end{array}\right],
\end{equation}
we obtain
\begin{equation}
W_k(M_k)=\tau\left(S-A\right)=\left[\begin{array}{cc}0&0\\-40.8241&-12.9855\end{array}\right],
\label{eq:wkm2}
\end{equation}
which also belongs to the range of the $k=-1$ branch.

As we move further to the left, the pairs of complex conjugate roots will always create matrices for which $W\left(m_{22}\right)\in\left(-\infty,\,-1\right)$. From the way in which $W_k(M)$ is created, we have that:
\begin{equation}
W_k\left(m_{22}\right)=\tau\left(2\Re\left(\lambda\right)-a_{22}\right).
\end{equation}
Since $a_{22}>0$ for this system and all the roots further to the left have negative real part with increasing absolute value, $W_k\left(m_{22}\right)$ will be always decreasing, thus remaining within the range of the $k=-1$ branch.
\section{Closing Remarks}
This short note presented a numerical example in which the full spectrum of a second order system can be found using only the $k=-1$ branch of the matrix Lambert w function, provided proper initial conditions. These results complement the remarks made by previous works on the applicability of the methodology under scrutiny.
The code used to create the example is available at \texttt{http://bit.ly/2uN8NZG} or can be requested via email to the author.

\end{document}